\documentclass[11pt]{amsart}

\usepackage{amsthm, amsfonts, amssymb,pinlabel,graphicx}
\usepackage[usenames,dvipsnames]{xcolor}
\usepackage{pinlabel}
\usepackage{algorithm}
\usepackage{algpseudocode}
\usepackage{bbm}
\usepackage{pgfplots}
\usepackage{longtable}
\usepackage[pdftex,%
  final,%
  colorlinks=true,%
  linkcolor=NavyBlue,%
  citecolor=NavyBlue,%
  filecolor=NavyBlue,%
  menucolor=NavyBlue,%
  urlcolor=NavyBlue,%
  bookmarks=true,%
  bookmarksdepth=3,%
  bookmarksnumbered=true,%
  bookmarksopen=true,%
  bookmarksopenlevel=2,%
]{hyperref}

\usetikzlibrary{patterns}

\numberwithin{equation}{section}
\graphicspath{{figs/}}

\theoremstyle{plain}
\newtheorem{thm}{Theorem}[section]
\newtheorem{cor}[thm]{Corollary}

\theoremstyle{definition}
\newtheorem{defn}[thm]{Definition}

\theoremstyle{remark}

\newtheorem{ex}[thm]{Example}

\newcommand{\zz}{\ensuremath{\mathbb{Z}}}
\newcommand{\qq}{\ensuremath{\mathbb{Q}}}

\newcommand{\dfn}[1]{{\textbf {#1}}}        % Definition

\DeclareMathOperator{\arf}{Arf}

\DeclareMathOperator{\rk}{rk}

\begin{document}

\title[Algorithmic Search for M\"obius Bands]{An Algorithmic Search for Knots Bounding M\"obius Bands in $B^4$} 

\begin{abstract}
This study examines the smooth non-orientable $4$-genus of knots using a large-scale computational search. Knot invariant obstructions were computed, and existing computational algorithms were extended to identify prime knots in $S^3$ of crossing number at most $14$ that bound M\"obius bands in $B^4$. The findings significantly expand the known examples of knots of non-orientable $4$-genus equal to $1$ and contribute new data relevant to the study of the non-orientable $4$-genus.
\end{abstract}

\date{\today}

\author[D.\ Lee]{Daniel Lee}
\address{Carnegie Mellon University} 
\email{dlee5@andrew.cmu.edu}

\author[J.M.\ Sabloff]{Joshua M. Sabloff} \address{Haverford College} \email{jsabloff@haverford.edu} 


\maketitle

% *******************
\section{Introduction}
\label{sec:intro}

While the study of slice knots and the orientable $4$-genus has long played an important role in low-dimensional topology, progress on its non-orientable analogue is of more recent vintage, with the first explicit definition not appearing until 2000 \cite{my:4-genus}.  Working in the smooth category, define the \dfn{non-orientable $4$-genus} $\gamma_4(K)$ of a knot $K \subset S^3$ to be the minimal first Betti number among all smooth, compact, properly embedded surfaces in $B^4$ with boundary $K$; if $\gamma_4(K)=1$, i.e.\ if one of those surfaces is a M\"obius band, we say that $K$ has a \dfn{M\"obius band filling}.

Early progress in exploring the non-orientable $4$-genus was made using lower bounds derived from classical knot invariants like the Arf invariant, the signature, and the linking form of the double branched cover \cite{gl:non-ori-4-genus, my:4-genus, viro:positioning, yasuhara:connecting}.  It was not until 2014, however, that the first example of a knot with $\gamma_4(K) > 3$ was found \cite{batson:non-ori-slice}.  Since then, new lower bounds using Heegaard-Floer theory \cite{batson:non-ori-slice, oss:unoriented}, Donaldson's diagonalization theorem \cite{jk:89}, Khovanov homology \cite{ballinger:concordance-kh}, and instanton homology \cite{ds:cs-clasp} have been developed, with more such research continuing to appear.

In this paper, we focus on searching for knots with M\"obius band fillings, adapting computer search techniques developed by Dunfield and Gong \cite{dg:slice} to find ribbon knots; see also \cite{gukov:mdp}.  Put another way, we investigate the output of algorithmic \emph{constructions}, though we also discuss how to make some of the \emph{obstructions} mentioned above readily computable. Overall, our goals are to develop computational tools and to provide raw data for future investigation. 

Prior efforts to survey the non-orientable $4$-genus of low-crossing knots include \cite{jk:89} for $8$ and $9$ crossing knots, \cite{ghanbarian:10} for $10$ crossing knots, and \cite{fairchild:11n} for non-alternating $11$ crossing knots; the constructions in these papers were, as far as we can tell, found by hand. 

In this study, we analyzed all $59,937$ prime knots of crossing number up $14$. Our algorithmic search yielded data that may be summarized as follows.

\begin{thm} \label{thm:overview}
    Of knots that have crossing number at most $14$, at least $24,504$  $(41\%)$ have M\"obius band fillings, while at least $22,412$ $(37\%)$ have $\gamma_4(K)>1$.  
\end{thm}

Of the knots determined to have M\"obius band fillings, $23,753$ are new examples. This collection includes $205$ knots with $11$ crossings (all alternating), $1,088$ with 12 crossings, $4,434$ with 13 crossings, and $18,026$ with 14 crossings.  See Section~\ref{sec:results} for a more detailed discussion. In all, the results described in this paper significantly expand the data set of the non-orientable $4$-genus of low-crossing knots.  

The paper is organized as follows. First, in Section~\ref{sec:background}, we review the Arf invariant and linking form obstructions, and present an algorithm for the using latter. Next, in Section~\ref{sec:band-algorithm}, we describe the computational model of Dunfield and Gong, along with our modifications for this study. In addition, we outline the band-specification conventions used throughout the paper and in our data set, which enable manual verification of our assertions of which knots have M\"obius band fillings.  Section~\ref{sec:results} contains a description of the output of the algorithm on knots up to crossing number $14$.  Finally, Appendix~\ref{sec:tables} lists knots with 11, 12, and 13 crossings with M\"obius band fillings which have yet to appear in the literature.

% **********
\subsection*{Acknowledgements}

We thank Sherry Gong, Fabian Ruehle, and Hailey Garcia for interesting and valuable discussions during the preparation of this paper.

% *******************
\section{Obstructions} 
\label{sec:background}

Our search for knots with M\"obius band fillings not only requires finding said fillings using the algorithm in Section~\ref{sec:band-algorithm}, below, but also eliminating from contention those knots that cannot have such a filling.  We use two readily computable obstructions: the Arf invariant obstruction of Yasuhara \cite{yasuhara:connecting} (see also \cite{gl:non-ori-4-genus}) and a computationally tractable version of the linking form obstruction of Murakami and Yasuhara \cite{my:4-genus}. Note that the Arf invariant provides an obstruction to smooth M\"obius bands, whereas the linking form obstructs topologically locally flat bands as well.

% **********
\subsection{The Arf Obstruction}
\label{ssec:arf}

The Arf invariant is a classical knot invariant based on a mod $2$ reduction of the Seifert form associated to a Seifert surface for $K$; for an exposition, see Lickorish \cite[Ch.\ 10]{lickorish:intro}. The Arf invariant of a knot can be determined by evaluating its Alexander polynomial at $-1$ \cite{murasugi:arf}:
\[\arf(K) = \begin{cases}
    0 & \Delta_K(-1) \equiv \pm 1 \mod 8 \\
    1 & \Delta_K(-1) \equiv \pm 3 \mod 8
\end{cases}\]

The obstruction arising from the Arf invariant is straightforward to state and to apply once the Arf invariant $\arf(K)$ and the signature $\sigma(K)$ of a knot $K$ have been computed.

\begin{thm}[Yasuhara \cite{yasuhara:connecting}, Gilmer-Livingston \cite{gl:non-ori-4-genus}] \label{thm:arf}
	If $\sigma(K) + 4\arf(K) \equiv 4 \mod 8$, then $K$ does not have a M\"obius band filling.
\end{thm}

\begin{ex}
	The Figure $8$ knot $K$ has $\sigma(K) = 0$ (since $K$ is amphichiral) and $\arf(K) = 1$.  Thus, since $\sigma(K) + 4 \arf(K) = 4$, the Figure $8$ knot does not have a M\"obius filling.  Since it is possible to construct a Klein bottle filling of $K$, we obtain $\gamma_4(K) = 2$. 
\end{ex}

% **********
\subsection{The Linking Form Obstruction}
\label{ssec:linking-obstruction}

We begin by following Murakami and Yasuhara \cite{my:4-genus} to define the linking form and set up the main obstruction. 

\begin{defn} \label{defn:linking-form}
Let $M$ be a $3$-manifold with finite first homology group.  The \dfn{linking form} $\lambda: H_1(M; \zz) \otimes H_1(M; \zz) \to \qq/\zz$
is defined by 
\[\lambda([x],[y]) = \frac{c \cdot y}{n},\]
where $nx = \partial c$.
\end{defn}

We say that a form is \dfn{metabolic} if there is a subgroup $H \subset G$ such that $|H|^2=|G|$ and $\lambda(h,h')=0$ for all $h,h' \in H$.

We now specialize to a knot $K \subset S^3$ and its double branched cover $D_K$.  It is a classical result that $\left|H_1(D_K;\zz)\right| = |\det(K)|$, so $D_K$ does have finite first homology.  Further, we may use $|\det(K)|$ as $n$ in the definition of the linking form.  That the linking form yields an obstruction to $K$ a M\"obius band filling comes from the following:

\begin{thm}[Corollary 2.7 of \cite{my:4-genus}] \label{thm:linking-obstruction}
	Let $K$ be a knot in $S^3$ that is either slice or bounds a M\"obius band in $B^4$.  Then the linking form $\lambda$ on $H_1(D_K;\zz)$ splits into a direct sum $(G_1, \lambda_1) \oplus (G_2, \lambda_2)$ where
\begin{enumerate}
\item $\lambda_1$ is represented by the $1 \times 1$ matrix $\bigl[\pm 1/|G_1|\bigr]$, and
\item $\lambda_2$ is metabolic.
\end{enumerate}
\end{thm}

The statement of the obstruction in Theorem~\ref{thm:linking-obstruction} requires some work in order to practically implement it.  In particular, we need to understand how to compute the linking form and how to test it for the existence of the required direct sum splitting.  For algorithmic ease, we choose to restrict to the following special case.

\begin{cor} 
\label{cor:linking-obstr}
Let $K$ be a knot in $S^3$ that is either slice or bounds a M\"obius band in $B^4$ whose determinant is square-free.  Then the linking form $\lambda$ is represented by the $1 \times 1$ matrix $\bigl[\pm 1/|G_1|\bigr]$ with respect to some generator of $H_1(D_K; \zz)$.
\end{cor}

To compute the linking form, we use the Goeritz matrix $U$ of a diagram of $K$.  The Goeritz matrix has the following two properties \cite{gl:signature}:
\begin{enumerate}
\item $U$ is a presentation matrix for $H_1(D_K;\zz)$, so $\det(K) = |\det(U)|$, and
\item $- U^{-1}$ defines $\lambda_K$. 
\end{enumerate}

\begin{figure}
    \begin{center}
        \includegraphics{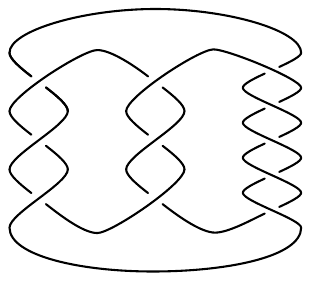}
    \end{center}
    \caption{The pretzel knot $P(3,3,-5)$.}
    \label{fig:P1}
\end{figure}

\begin{ex} \label{ex:P-linking}
	The pretzel knot $P = P(3,3,-5)$ pictured in Figure~\ref{fig:P1} has Goeritz matrix 
	\[U = \begin{bmatrix} -6 & 3 \\ 3 & 2 \end{bmatrix}.\]
	We compute that $|\det U| = 21$, and that 
	\[-U^{-1} = \begin{bmatrix}\frac{2}{21} & -\frac{1}{7} \\ -\frac{1}{7} & -\frac{2}{7} \end{bmatrix}.\]
	Since $\det U$ is square-free, we see that if $K$ were to bound a M\"obius band, we would need to present the linking form $\lambda$ of $G$ as a $1 \times 1$ matrix of the form $\bigl[\pm 1/21\bigr]$. Thus, we need to find all possible generators $g$ of $G$ and test the linking form to see if we obtain $\lambda(g,g) = \pm \frac{1}{21}$.  Of course, we really need only find a single generator, and then look at its multiples to test all possible generators.
	
    We claim that, with respect to the presentation of $G$ by $U$, $g_0 = \left[\begin{smallmatrix} 1 \\ 0 \end{smallmatrix}\right]$ is a generator for $G$. To prove the claim, we check that none of $g$, $3g$, or $7g$ vanish in $G$ by ensuring that $U x = q g_0$  does not have integral solutions for $q\in \{1,3,7\}$.
	
	To test the generators, we compute that $\lambda(ng_0,ng_0) = \frac{2n^2}{21}$.  A straightforward computation shows that $2n^2 \not\equiv \pm 1 \mod 21$ for $n \in \{1, \ldots, 20\}$. Hence the knot $P$ is obstructed from having a M\"obius band filling by Corollary~\ref{cor:linking-obstr}.
\end{ex} 

\begin{algorithm}[ht!] 
\caption{Obstructing M\"obius Fillings Using Linking Forms}
\label{alg:linking}
\begin{algorithmic}[1]

\Require Goeritz matrix $U$ of $K$
\Ensure Whether $K$ is obstructed from having a M\"obius band filling by the linking form

\State $V \gets -U^{-1}$
\State $n \gets \rk U$
\State $\Delta \gets |\det U|$

\Statex

\Function{FindGenerator}{$m$,$A$}
	\State $r_1, \cdots, r_k \gets$ proper factors of $|\det A|$
	\For{$i$ from $1$ to $m$}
		\If{Solution of $Ax = r_je_i$ is integral for some $j \in \{1, \ldots, k\}$}
			\State \textbf{break}
		\Else
			\State \Return i \Comment{$e_i$ is a generator}
		\EndIf
	\EndFor
	\State \Return None \Comment{No $e_i$ is a generator} 
\EndFunction

\Statex

\If{Prime factorization of $\Delta$ has repeated factors}
	\Return ``$\det U$ not square-free''
\Else
	\State $i \gets$ FindGenerator($n$,$U$)
    \If{$i=$ None} \Return ``Need to find generator by hand''
    \EndIf
	\State $\ell \gets V_{ii} \Delta$
	\For{$j$ from $1$ to $\Delta-1$}
		\If{$j^2 i \equiv \pm 1 \mod \Delta$}
			\State \Return ``M\"obius band Possible''
		\EndIf
	\EndFor
	\State \Return ``M\"obius band Obstructed''
\EndIf
\end{algorithmic}
\end{algorithm} 

We generalize the process of Example~\ref{ex:P-linking} to a computational algorithm that we may use to filter out knots that do not have M\"obius fillings. In Algorithm \ref{alg:linking}, we denote by $e_i$ the $i^{th}$ basis vector of $\zz^n$.  There are two caveats to the application of this algorithm:  first, $\det K$ must be square-free.  Second, one of the $e_i$ must be a generator of $H_1(D_K; \zz)$, as presented by the Goeritz matrix, which may not always be the case.

% *******************
\section{Band-Finding Algorithm}
\label{sec:band-algorithm}

In \cite{dg:slice}, Dunfield and Gong introduced a computational framework for detecting slice knots through iterative searches over band moves. Briefly, their algorithm attempts to find a slicing disk for a given knot by systematically enumerating a collection of diagrammatically defined oriented bands, ordered from the simplest to the most complicated. For each candidate band, the algorithm performs the corresponding band move and examines whether the resulting link is an unlink comprised of previously detected ribbon knots. The algorithm terminates if such a band is found, or else it continues the search with increasingly complex bands up to a complexity threshold. We adapt Dunfield and Gong's algorithm to the non-orientable setting by using unoriented, rather than oriented, bands.

This section elaborates on the summary above and presents a modification that searches for unoriented bands that yield slice knots when attached to a given knot. Along with the description, we provide a running example with the $5_2$ knot. For implementation, we utilize SnapPy \cite{SnapPy} and SageMath \cite{sagemath}. The code is available at
\begin{center} \url{https://github.com/LeeJoPing/Mobius-Band-Search}
\end{center}

% **********
\subsection{Finding Candidate Bands}
\label{ssec:band-search}

\begin{figure}
\centering
\begin{minipage}{0.4\textwidth}
  \centering
  \includegraphics[width=\linewidth]{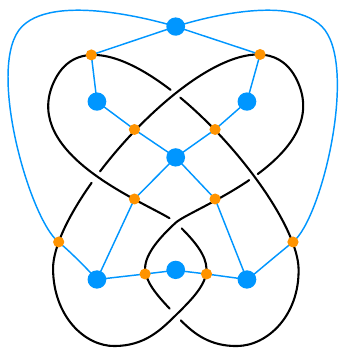}
  \small (a)
\end{minipage}
\hspace{0.1\textwidth}
\begin{minipage}{0.4\textwidth}
  \centering
  \includegraphics[width=\linewidth]{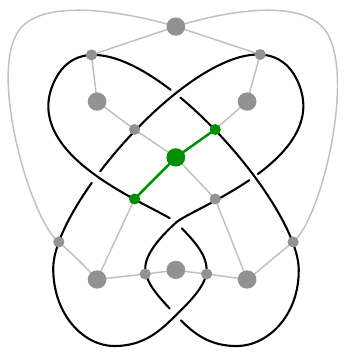}
  \small (b)
\end{minipage}
\caption{(a) A modified dual graph representation of a three-twist ($5_2$) knot.
(b) A minimal simple path (shown in green) in $G_D$ that can be used to attach an unoriented band.}
\label{fig:52knot}
\end{figure}

To make the ribbon search computationally tractable, Dunfield and Gong \cite{dg:slice} restrict the search space to a finite, diagrammatically defined collection of candidate bands. Given a knot diagram, the algorithm begins by constructing the dual graph $G_D$: each vertex of $G_D$ corresponds to a distinct region of the planar projection, and each edge connecting two vertices represents the arc shared by the paired regions. To encode the crossing information required for reconstructing band attachments, $G_D$ is modified by adding a vertex in the middle of each edge where it intersects the knot diagram. Figure~\ref{fig:52knot}(a) illustrates an example of such a modified dual graph. 

\begin{figure}
\labellist
\small\hair 2pt
 \pinlabel {$0$} [b] at 16 35
 \pinlabel {$1$} [l] at 33 18
 \pinlabel {$2$} [b] at 47 35
 \pinlabel {$3$} [l] at 33 54
 \pinlabel {$0$} [b] at 127 35
 \pinlabel {$1$} [l] at 143 13
 \pinlabel {$2$} [b] at 155 35
 \pinlabel {$3$} [l] at 143 49
\endlabellist
    \centering
    \includegraphics{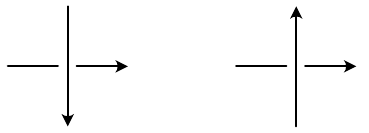}
    \caption{The conventions for strand index at a crossing depend on the sign of the crossing.}
    \label{fig:strandidx}
\end{figure}
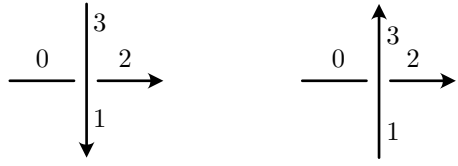

Each arc in the knot diagram is encoded as a pair \texttt{(c, s)}, where \texttt{c} denotes the crossing label and \texttt{s} indicates the strand index at that crossing. The strand indices at each crossing follow the convention inherited from SnapPy; see Figure~\ref{fig:strandidx}. The crossing labels are also inherited from SnapPy's database. A band specification starts as a simple path in $G_D$ with endpoints on the knot diagram. This path is encoded as a tuple consisting of the starting arc, the sequence of arcs traversed, and the ending arc. To reduce computational complexity and time, the search algorithm prioritizes bands of shorter length, which are determined by the number of edges a band traverses in a dual graph, and limits the search by imposing a maximum allowable length. Figure~\ref{fig:52knot}(b) shows a minimal-length path corresponding to a candidate band for the knot $5_2$.

Once a path is fixed, the algorithm enumerates the allowable twist parameters. Let $tw \in \zz$ denote the twisting parameter, where a positive $tw$ represents $tw$ clockwise twists and a negative value represents $|tw|$ counterclockwise twists. Let $p$ denote the relative parity of the orientations of the arcs at either end of a band:  $p=1$ if the orientations of the attaching arcs match and $p=0$ otherwise. In the oriented setting described in \cite{dg:slice}, the number of twists is constrained by a parity condition of the form $tw - p \equiv  0 \pmod 2$.  To adapt the algorithm for unoriented band attachments, the parity condition is modified to $tw - p \equiv 1 \pmod 2$. In our implementation, the twist parameter is limited to at most three half-twists in either direction to constrain the search. 

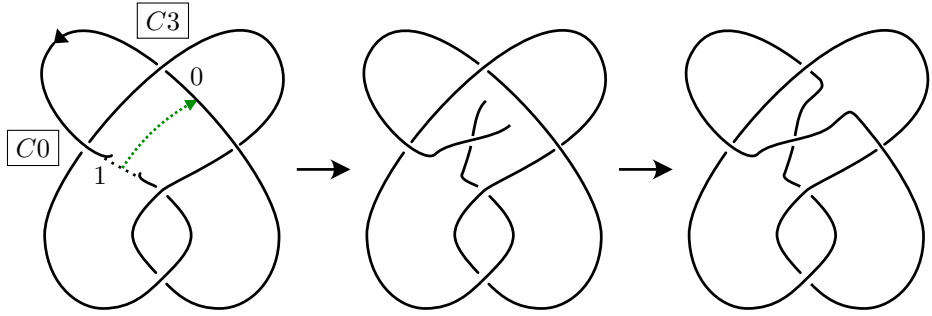
\begin{figure}
\labellist
\small\hair 2pt
 \pinlabel {\fbox{$C0$}} [r] at 12 63
 \pinlabel {$1$} [t] at 25 58
 \pinlabel {\fbox{$C3$}} [b] at 48 101
 \pinlabel {$0$} [l] at 56 89
\endlabellist
    \centering
    \includegraphics[width=0.95\linewidth]{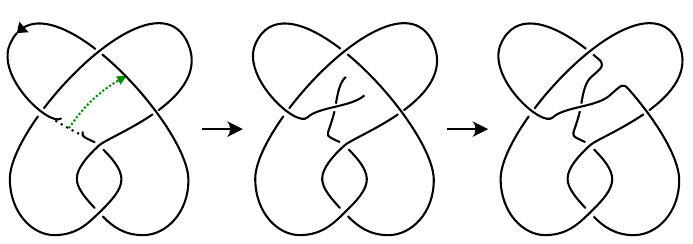}
    \caption{The attachment of a (twisted) unoriented band to the $5_2$ knot. The crossing and strand numbers are labeled on the leftmost figure}
    \label{fig:52bandaction}
\end{figure}   

Figure \ref{fig:52bandaction} depicts the process of attaching an unoriented band to a diagram of the $5_2$ knot. In this example, the algorithm applies a single clockwise twist to break the orientation.

\begin{figure}
\labellist
\small\hair 2pt
 \pinlabel {\fbox{$C0$}} [tl] at 135 61
 \pinlabel {$1$} [l] at 133 75
 \pinlabel {\fbox{$C3$}} [tl] at 97 61
 \pinlabel {$2$} [l] at 95 75
 \pinlabel {$3$} [t] at 77 61
\endlabellist
    \centering
    \includegraphics{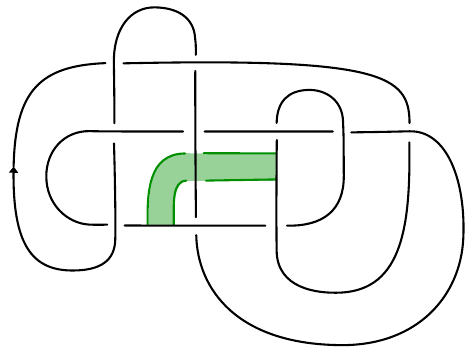}
    \caption{A band that yields the $6_1$ knot (which is slice) when attached to the knot $10a31$.}
    \label{fig:bandaction}
\end{figure}
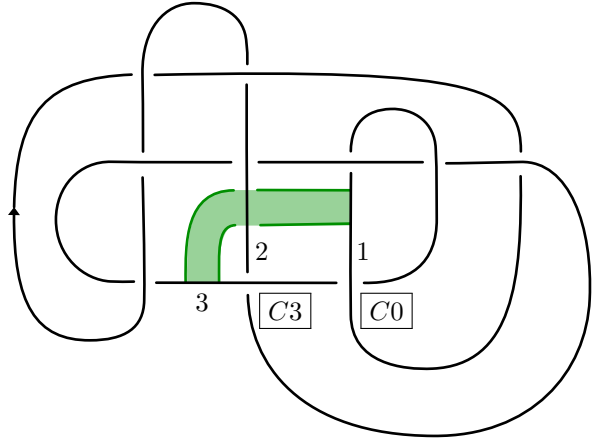

The final step in enumerating candidate bands involves specifying over and under information for each intermediate arc in the knot diagram. At each intersection between a band and an arc, the band may pass over or under the arc, resulting in $2^n$ possible configurations when the path intersects $n$ intermediate arcs. The algorithm iterates through bands with all combinations of over and under-passing along the path. Figure ~\ref{fig:bandaction} presents an example of a band passing under an arc during a band move.

Within this framework, a band is defined as a three-component data structure that captures all the information discussed above. The first component represents the band path. The second component is a list of Boolean values indicating whether the band passes above each intermediate arc (\texttt{True} indicates that the band passes above the respective arc, and \texttt{False} denotes that it lies below). If the band does not pass through any intermediate arcs, this list is empty. The third component specifies the twisting parameter, indicating the number of clockwise or counterclockwise twists applied to the band. 

\begin{ex}
The band that unknots the $5_2$ knot in Figure~\ref{fig:52bandaction} is expressed as
\begin{center}
  \texttt{Band([(0, 1), (3, 0)], [], 1)}.
\end{center}
This means that the band starts at the arc with crossing label $0$ and strand index $1$ and does not pass through any intermediate arcs. Before attachment, the band is twisted clockwise once, then ends at the arc with the crossing label $3$ and strand index $0$. 
\end{ex}

\begin{ex}
The band pictured in Figure~\ref{fig:bandaction} is denoted by 
\begin{center}
    \texttt{Band([(0,1), (3,2), (3,3)], [False], 0)}.
\end{center}
\end{ex}

We provide an overview of the procedure for enumerating all candidate bands in Algorithm \ref{alg:UnorientedBands}. Note that the band attachment implementation is based on the algorithm used in \cite{dg:slice}.

\begin{algorithm}
\caption{Minimal Length Unoriented Bands}
\label{alg:UnorientedBands}
\begin{algorithmic}[1]

\Require A link diagram $L$ and maximum length parameter $\text{max\_band\_len}$
\Ensure A collection $\mathcal{B}$ of allowable minimal unoriented bands

\State Construct the modified dual graph $G_D$ of $L$
\State Compute shortest-path distances for all-pairs in $G_D$
\State $\mathcal{B} \gets \emptyset$

\For{each component $C$ of $L$}
    \For{each unordered pair of arcs $(a_1,a_2)$ in $C$}

        \State $d \gets$ shortest-path distance in $G_D$ between $a_1$ and $a_2$
        \If{$d > \text{max\_band\_len}$}
            \State \textbf{continue}
        \EndIf

        \For{each minimal path $\gamma $ in $ G_D$ connecting $a_1$ to $a_2$}
            \State Construct a band core on $\gamma$
            \State $\mathcal{T} \gets $ set of twist numbers that break orientation
            
            \For{each twist number $t \in \mathcal{T}$}
                \For{each over/under passing}
                    \State $B \gets$ band defined by $\gamma$ , $t$, and crossing data
                    \State $\mathcal{B} \gets \mathcal{B} \cup \{B\}$
                \EndFor
            \EndFor

        \EndFor

    \EndFor
\EndFor

\State Sort $\mathcal{B}$ by increasing path length
\State \Return $\mathcal{B}$

\end{algorithmic}
\end{algorithm}

% **********
\subsection{Testing Bands}
\label{ssec:test-band}

Once a band is fully specified and attached, the original Dunfield and Gong algorithm analyzes the resulting link to determine whether it is an unlink of unknots or of known ribbon knots \cite{dg:slice}. The search process terminates if either outcome is achieved. If not, the algorithm proceeds with the remaining candidate bands until such candidates have been exhausted. In our unoriented setting, the modification guarantees that the resulting object remains a one-component link. Thus, each knot resulting from an unoriented band move is evaluated using the same criteria to determine whether it is trivial or ribbon. If an unoriented band transforms the given knot into a slice knot, it is concluded that the knot bounds a M\"obius band, and the band is returned. 

% *******************
\section{Results}
\label{sec:results}

% **********
\subsection{Inputs}
\label{ssec:input}

We applied our algorithm to knots with crossing numbers from $3$ to $14$. The non-orientable 4-genus for knots with up to 10 crossings and for non-alternating knots with 11 crossings has been determined in previous work \cite{fairchild:11n,ghanbarian:10,jk:89}. We use these data to validate the model. For alternating knots with 11 crossings and all knots with 12 to 14 crossings, we investigate previously unexplored examples of knots that bound a M\"obius band.

For knots with crossing numbers from $3$ to $13$, we collected data from KnotInfo \cite{knotinfo:data} and their diagrams from SnapPy \cite{SnapPy}. For knots with $14$ crossings, we used SnapPy \cite{SnapPy} to obtain both the knot diagrams and their invariants, and then computed the Arf invariants using the method described in Section~\ref{ssec:arf}. Altogether, the dataset includes $59,\!937$ knots. 

Our implementation follows the Dowker-Thistlethwaite (DT) naming convention used by KnotInfo \cite{knotinfo:data}. Specifically, knots are indexed according to the lexicographic ordering of their minimal DT notations, with letter `a' and `n' denoting alternating and non-alternating knots, respectively. For example, the trefoil knot $3_1$ is denoted as ``K3a1''. The knots are assumed to follow the standard orientation convention in which the first crossing is right-handed. 

% **********
\subsection{Obstructions}
\label{ssec:obstructions}

The first step was to filter out knots that are obstructed from having a M\"obius band filling using the Arf invariant (Theorem \ref{thm:arf}) and the linking form (Theorem \ref{thm:linking-obstruction}). Together, these obstructions filtered $22,412$ of the $59,937$ knots of at most $14$ crossings. See Table~\ref{tab:obstruction} for more details. We note that the two invariants consistently obstructed about a third of all knots at each level of crossing number.

\begin{table}[ht]
\centering
\begin{tabular}{c|r||r|r|r}
 \hline
 Crossing \# & \# Knots & Arf & Linking Form  & Total Obstructed \\
 \hline
3-8 & 35 & 11 & 11 & 12 \\ 
 9 & 49 & 12 & 13 & 16 \\ 
 10 & 165 & 43 & 54 & 59 \\ 
 11 & 552 & 141 & 165 & 189 \\  
 12 & 2176 & 545 & 707 & 807 \\ 
 13 & 9988 & 2417 & 3190 & 3630 \\ 
 14 & 46972 & 11415 & 15431 & 17699 \\ \hline
 Total & 59937 & 14584 & 19571 & 22412 \\
\end{tabular}

\vspace{.1in}

\caption{Number of knots of each crossing number together with the numbers obstructed by the Arf invariant (Theorem \ref{thm:arf}), by the linking form (Theorem \ref{thm:linking-obstruction}), and the total number  obstructed by either invariant.}
\label{tab:obstruction}
\end{table}

% **********
\subsection{Algorithm Output}
\label{ssec:output}

Out of $59,937$ knots of at most $14$ crossings, $37,525$ were not obstructed from bounding a M\"obius band and hence were tested with the band-finding algorithm.  We ran the algorithm twice:  first comparing the knots resulting from a band attachment with the existing ribbon knot database, and the second time using the algorithm of \cite{dg:slice} to check for sliceness.  In the second run, we only checked the following knots:
\begin{itemize}
\item 11a knots that remained unknown after the first run,
\item and all known slice knots with at most $13$ crossings but for which the algorithm did not find a band, and
\item those knots reported in previous literature as having $\gamma_4 = 1$ but for which we had not found a band.
\end{itemize} 
We did not run the more thorough check on the full set of knots due to the likely significant increase in computational time.

The algorithm found $24,504$ knots with at most $14$ crossings with M\"obius band fillings; the vast majority of these knots were unknown to this point. The algorithm (and the obstructions) were unable to determine whether $13,021$ knots bound a M\"obius band.  As before, the underlying data that lists the knots that bound M\"obius bands, as well as the bands themselves described using the notation set out in Section~\ref{ssec:band-search}, can be found at:
\begin{center} \url{https://github.com/LeeJoPing/Mobius-Band-Search}
\end{center}

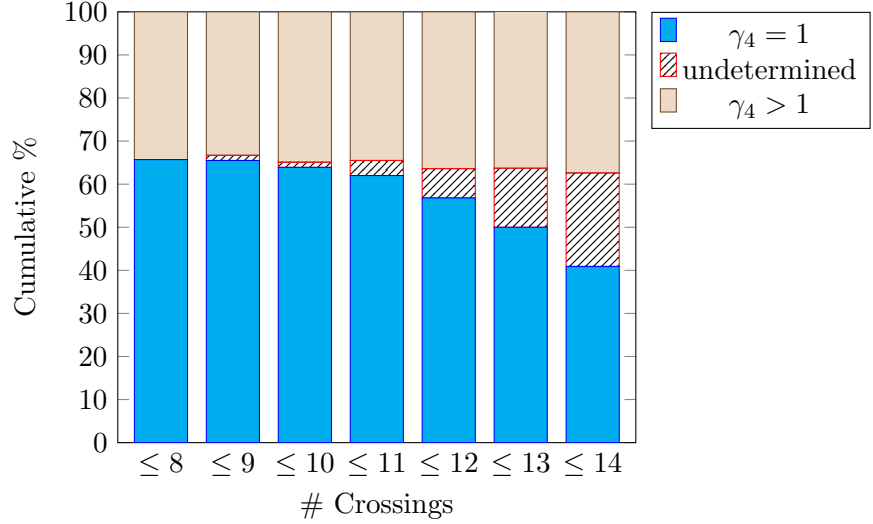
\begin{figure}
\begin{tikzpicture}
\begin{axis}[
ybar stacked,
bar width = 20,
ytick distance=10,
legend pos = outer north east,
ylabel={Cumulative \%},
symbolic x coords={$\leq 8$, $\leq 9$, $\leq 10$, $\leq 11$, $\leq 12$, $\leq 13$, $\leq 14$},
xtick=data,
xlabel={\# Crossings},
ymin={0},
ymax={100}
%x tick label style={rotate=45,anchor=east},
]
\addplot+ [ybar,fill=cyan] coordinates {($\leq 8$,65.7) ($\leq 9$,65.5)
($\leq 10$,63.9) ($\leq 11$,62) ($\leq 12$,56.8) ($\leq 13$,50) ($\leq 14$,40.9)};
\addplot+ [ybar,pattern=north east lines] coordinates {($\leq 8$,0) ($\leq 9$,1.2)
($\leq 10$,1.2) ($\leq 11$,3.5) ($\leq 12$,6.8) ($\leq 13$,13.7) ($\leq 14$,21.7)};
\addplot+ [ybar] coordinates {($\leq 8$,34.3) ($\leq 9$,33.3)
($\leq 10$,34.9) ($\leq 11$,34.5) ($\leq 12$,36.4) ($\leq 13$,36.4) ($\leq 14$,37.4)};
\legend{$\gamma_4 = 1$, undetermined, $\gamma_4>1$}
\end{axis}
\end{tikzpicture}
\caption{Cumulative percentage of knots according to their non-orientable $4$-genus.}
\label{fig:cumu-mobius}
\end{figure}

\begin{figure}
\begin{tikzpicture}
\begin{axis}[
%ybar,
bar width = 20,
ymajorgrids=true,
legend pos = outer north east,
ylabel={Cumulative \%},
symbolic x coords={$\leq 8$, $\leq 9$, $\leq 10$, $\leq 11$, $\leq 12$, $\leq 13$},
xtick=data,
ytick distance=10,
xlabel={\# Crossings},
ymin={0},
ymax={80}
%x tick label style={rotate=45,anchor=east},
]
\addplot [ybar, fill=cyan] coordinates {($\leq 8$,65.7) ($\leq 9$,65.5)
($\leq 10$,63.9) ($\leq 11$,62) ($\leq 12$,56.8) ($\leq 13$,50)};
\addplot+ [sharp plot, red] coordinates {($\leq 8$,11.4) ($\leq 9$,8.3)
($\leq 10$,8.4) ($\leq 11$,6.4) ($\leq 12$,5.3) ($\leq 13$,3.9)};
\addplot+ [sharp plot, red] coordinates {($\leq 8$,65.7) ($\leq 9$,61.9)
($\leq 10$,63.9) ($\leq 11$,62.8) ($\leq 12$,63.3) ($\leq 13$,62.7)};
\legend{$\gamma_4 = 1$, $g_4=0$,$g_4=1$}
\end{axis}
\end{tikzpicture}
\caption{Cumulative percentage of knots to compare their non-orientable and orientable $4$-genera.  Orientable genera are exact numbers, while the non-orientable genus indicates known knots that bound M\"obius bands.}
\label{fig:comp-genus}
\end{figure}
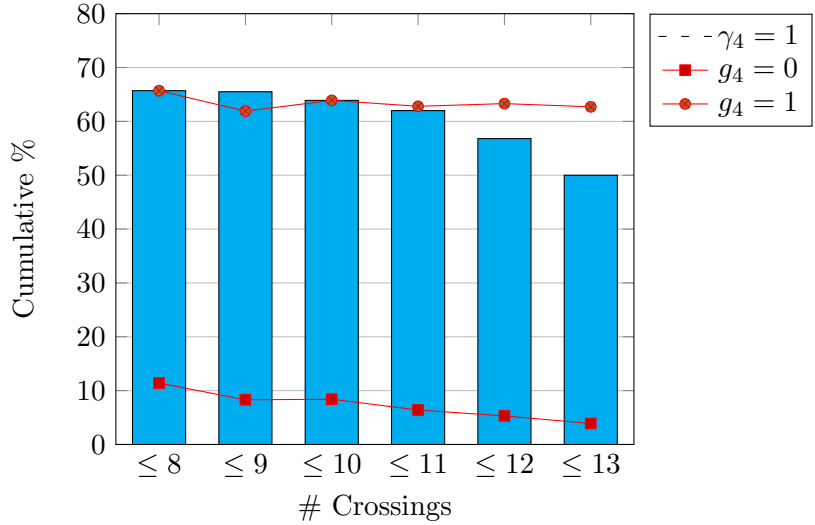

For a first look at the data, we turn to Table~\ref{tab:mobius}, which lists the numbers of knots that do or do not bound M\"obius bands, or for whose status is undetermined listed by crossing number.  Figure~\ref{fig:cumu-mobius} depicts the cumulative percentages of these three classes of knots. Finally, Figure~\ref{fig:comp-genus} compares the cumulative percentage of knots that bound M\"obius bands to the cumulative percentages of those that have (orientable) $4$-genus $0$ and $1$. 

\begin{table}
\centering
\begin{tabular}{c|r||r|r|r}
\hline
Crossing \# & \# Knots & $\gamma_4 = 1$ & $\gamma_4 > 1$ & Undetermined \\
\hline
3-8 & 35 & 23 & 12 & 0 \\ 
9 & 49& 32 & 16 & 1 \\ 
10 & 165& 104 & 59  & 2 \\ 
11 & 552& 338 & 189 & 25 \\ 
12 & 2176& 1194 & 807 & 175 \\ 
13 & 9988& 4787 & 3630 & 1571 \\ 
14  & 46972& 18026 & 17699 & 11247 \\ \hline
Total & 59937 & 24504 & 22412 &13021 
\end{tabular}
\caption{Number of knots at each crossing number according to their non-orientable $4$-genus, $\gamma_4$, after our computation.}
\label{tab:mobius}
\end{table}

% ***
\subsubsection{9 and 10 Crossing Knots}

The algorithm was unable to determine whether three knots of crossing number at most 10 bound M\"obius bands: K9a37 ($9_{40}$ in the standard tables), K10a59 ($10_2$), and K10a81 ($10_{46}$).  In each case, one may use an argument based on Donaldson's diagonalization theorem for definite intersection forms to determine that none of these knots bound a M\"obius band.  For details of this argument, see \cite{jk:89} for K9a37 and \cite{ghanbarian:10} for the 10-crossing knots.

% ***
\subsubsection{11n Knots}

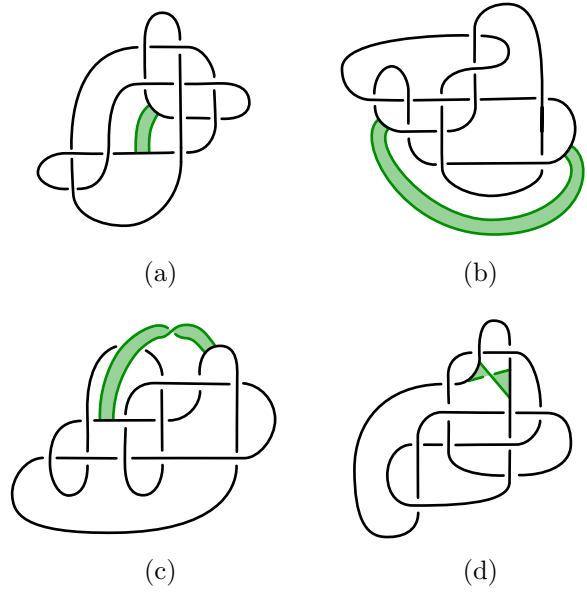
\begin{figure}
\labellist
\small\hair 2pt
\pinlabel {(a)} [b] at 60 112
\pinlabel {(b)} [b] at 180 112
\pinlabel {(c)} [b] at 60 0
\pinlabel {(d)} [b] at 180 0
\endlabellist
    \centering
    \includegraphics{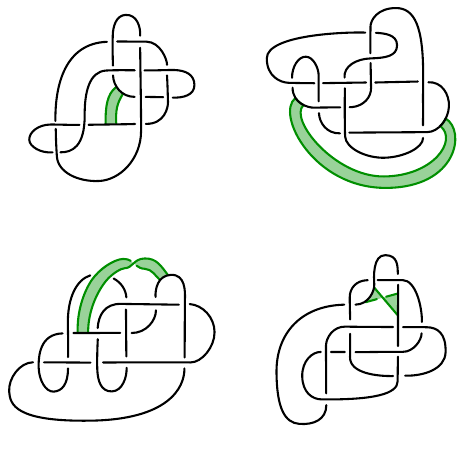}
    \caption{Unoriented bands to slice knots for (a) 11n92 (which yields the unknot), (b) 11n125 (which yields 8n1), (c) 11n131 (which yields 6a3), and (d) 11n161 (which yields 6a3).}
    \label{fig:11n-bands}
\end{figure}

The algorithm yielded several corrections to the calculations of the non-orientable $4$-genus in \cite{fairchild:11n} (and hence in \cite{knotinfo:data}), though it validated the vast majority of the computations in \cite{fairchild:11n}.  On one hand, the knots 11n92, 11n125, 11n131, and 11n161 all bound M\"obius bands; see Figure~\ref{fig:11n-bands}.

In the other direction, the algorithm failed to find M\"obius bands for 11n36 and 11n57, which were claimed to have $\gamma_4=1$ in \cite{fairchild:11n}.  The band drawn in  \cite{fairchild:11n} for 11n36  is claimed to produce the knot $10_{129}$.  We computed, however, that the resulting knot has a different Alexander polynomial than $10_{129}$ and, in fact, is the 11n82.  The signature of 11n82 is nonzero, hence it is not a slice knot. Similarly, the band drawn in \cite{fairchild:11n} for 11n57 is claimed to yield the unknot. Our computations, however, show that the result is the figure 8 knot, which is not slice.

That said, we were not able to definitively obstruct either the 11n36 or the 11n57 from bounding a M\"obius band using the Arf invariant, the linking form, or an argument using Donaldson's diagonalization theorem.  Thus, we list these knots as having undetermined $\gamma_4$.

% ***
\subsubsection{Future Directions}

These computations suggest several promising directions for future investigation. For example, one may ask how the proportion of knots obstructed from bounding an M\"obius band or the cumulative proportion of knots with $\gamma_4 = 1$ change as the crossing number increases. Additionally, comparing these trends with the distribution of knots with $g_4 = 1$ may yield further insights. Although several potential trends may be posited from these tables and figures, one should be exceedingly careful not to attempt to formulate general laws from observations of low-crossing knots alone. 

Addressing these questions will require extending computations to knots with substantially higher crossing numbers, making improvements to the computational model increasingly important. One potential approach is to incorporate machine learning techniques to guide the search through the vast solution space more efficiently. In \cite{gukov:mdp}, a Bayesian-optimized Markov decision process is implemented to identify ribbon disks in arbitrary knots. Adapting a similar framework to the non-orientable setting could significantly enhance the efficiency of searching for M\"obius band constructions. In addition, parallelization through multithreading or other high-performance computing methods could substantially reduce computation times, enabling the model to scale more effectively to larger knot families and higher crossing numbers. Finally, we suspect that the growing number of undetermined knots owes more to the relative inadequacy of the obstructions as opposed to a failure to find useful band moves.  Thus, sharpening the linking form algorithm or devising diagrammatic computational methods for other obstructions, such as the $\upsilon$ and $d$ invariants, would also yield a significant impact.

\appendix 

% ********************
\section{Tables of Knots Bounding M\"obius Bands with 11 (alternating), 12, and 13 Crossings}
\label{sec:tables}

The tables below list new examples of knots with 11 (alternating only), 12, or 13 crossings that have M\"obius band fillings.  There are $6,195$ of these.  See the Github repository mentioned in Section~\ref{ssec:output} for a list of $18,026$ knots of crossing number $14$ with M\"obius band fillings, as well as for a list of knots with crossing number at most $14$ with $\gamma_4>1$.

% **********
\subsection{11a}

\begin{center}
\begin{longtable}{|c|c|c|c|c|c|c|c|c|c|c|c|c|}
\hline
1&3&4&7&10&11&12&15&16&18&19&21&22 \\
\hline
23&25&27&28&31&32&33&34&35&36&42&43&45 \\
\hline
50&52&54&55&57&58&60&61&62&63&65&68&70 \\
\hline
72&73&75&76&77&78&80&81&85&87&90&91&92 \\
\hline
94&96&97&99&100&101&102&103&105&107&108&110&111 \\
\hline
113&114&115&117&118&120&121&122&123&124&125&126&128 \\
\hline
130&131&132&133&134&135&136&137&138&139&141&142&145 \\
\hline
146&147&148&149&150&151&153&154&155&156&157&158&159 \\
\hline
162&164&165&166&167&169&170&172&174&175&180&181&183 \\
\hline
186&191&193&194&196&199&201&202&204&205&207&210&211 \\
\hline
215&216&217&218&219&221&223&225&226&227&229&230&231 \\
\hline
232&235&236&237&239&242&243&244&245&247&250&252&254 \\
\hline
257&258&260&262&263&266&267&268&269&270&271&281&282 \\
\hline
283&284&285&286&291&294&295&296&298&300&301&304&307 \\
\hline
309&310&311&312&316&320&322&323&324&326&327&328&333 \\
\hline
334&335&336&340&341&342&343&346&347&353&356&357&358 \\
\hline
359&361&363&364&365&367 & -- & -- &-- &-- &-- &-- &-- \\
\hline
\end{longtable}
\end{center}

% **********
\subsection{12a}

\begin{center}
\begin{longtable}{|c|c|c|c|c|c|c|c|c|c|c|}
\hline
2&3&5&7&10&11&12&16&17&18&19 \\
\hline
20&21&27&28&31&32&34&37&40&42&45 \\
\hline
46&50&53&54&57&59&61&62&63&65&69 \\
\hline
71&73&74&75&77&78&79&81&82&83&84 \\
\hline
86&87&89&90&91&92&93&94&95&98&99 \\
\hline
100&102&104&105&107&110&111&112&118&119&122 \\
\hline
123&124&127&129&130&131&135&138&139&141&144 \\
\hline
145&147&149&150&153&155&158&159&160&161&164 \\
\hline
166&168&169&170&173&174&175&177&178&179&182 \\
\hline
183&184&187&189&192&194&196&198&199&203&205 \\
\hline
210&211&213&214&215&216&219&220&221&222&225 \\
\hline
226&227&228&229&230&231&233&235&236&238&240 \\
\hline
241&242&244&245&247&249&250&251&252&255&256 \\
\hline
257&258&260&261&262&263&265&266&267&270&271 \\
\hline
272&275&278&279&281&282&283&284&285&286&287 \\
\hline
290&291&292&293&294&295&298&303&306&308&311 \\
\hline
312&315&316&317&318&319&320&321&324&325&326 \\
\hline
328&329&333&334&336&337&339&343&344&348&349 \\
\hline
350&353&355&356&360&361&362&363&364&365&366 \\
\hline
368&369&375&376&377&378&379&382&383&384&385 \\
\hline
386&391&393&394&395&399&401&402&403&404&405 \\
\hline
406&407&408&412&413&415&416&417&419&420&421 \\
\hline
422&423&424&425&427&429&431&432&433&434&435 \\
\hline
436&438&439&440&441&442&443&445&447&448&450 \\
\hline
451&454&455&456&458&459&464&465&468&469&472 \\
\hline
473&475&477&478&479&481&484&486&488&489&490 \\
\hline
494&495&496&497&498&499&500&501&509&510&511 \\
\hline
512&513&515&519&520&522&524&526&527&530&531 \\
\hline
534&535&537&538&540&543&548&549&552&555&557 \\
\hline
559&562&563&564&566&567&568&572&576&577&579 \\
\hline
580&581&582&583&586&591&594&596&598&600&601 \\
\hline
604&606&607&608&609&615&617&619&620&621&623 \\
\hline
624&625&626&629&630&631&634&635&638&642&643 \\
\hline
644&646&648&649&652&658&659&660&662&663&666 \\
\hline
667&671&672&676&677&678&680&682&683&684&687 \\
\hline
690&695&704&710&711&713&715&718&721&723&724 \\
\hline
725&727&729&734&735&739&741&743&746&747&748 \\
\hline
749&751&752&754&755&756&757&758&759&760&761 \\
\hline
764&768&770&771&774&775&776&779&782&784&786 \\
\hline
790&791&794&797&798&799&800&814&819&821&823 \\
\hline
824&825&826&833&835&838&839&842&843&845&847 \\
\hline
849&851&852&856&857&860&863&866&867&870&871 \\
\hline
875&877&879&883&886&887&889&892&894&896&897 \\
\hline
898&900&901&908&910&911&913&914&915&917&918 \\
\hline
920&921&922&924&926&929&930&934&936&938&939 \\
\hline
941&945&947&949&951&952&954&955&958&959&960 \\
\hline
962&964&965&966&967&968&971&972&975&977&978 \\
\hline
979&983&985&988&989&992&993&995&999&1000&1005 \\
\hline
1007&1010&1011&1016&1017&1018&1019&1020&1023&1027&1028 \\
\hline
1029&1032&1033&1034&1038&1039&1040&1047&1053&1054&1055 \\
\hline
1057&1061&1062&1064&1066&1067&1068&1071&1073&1077&1081 \\
\hline
1083&1085&1086&1087&1089&1090&1093&1095&1096&1098&1101 \\
\hline
1104&1105&1108&1110&1111&1113&1118&1119&1125&1126&1127 \\
\hline
1128&1132&1133&1134&1135&1136&1138&1140&1141&1144&1145 \\
\hline
1147&1148&1149&1151&1154&1158&1161&1162&1163&1165&1166 \\
\hline
1168&1169&1170&1174&1175&1176&1178&1181&1182&1185&1186 \\
\hline
1188&1190&1194&1195&1196&1198&1202&1203&1205&1207&1211 \\
\hline
1215&1216&1219&1220&1221&1224&1225&1226&1227&1230&1231 \\
\hline
1233&1234&1236&1237&1238&1239&1241&1243&1247&1249&1250 \\
\hline
1253&1257&1262&1264&1267&1269&1271&1272&1277&1280&1282 \\
\hline
1283&1284&1286& -- & -- &-- &-- &-- &-- &--&-- \\
\hline
\end{longtable}
\end{center}

% **********
\subsection{12n}

\begin{center}
\begin{longtable}{|c|c|c|c|c|c|c|c|c|c|c|c|c|}
\hline
1&2&3&4&6&8&9&10&12&14&15&16&17 \\
\hline
19&20&21&23&24&26&28&29&31&34&38&40&42 \\
\hline
43&44&46&48&49&50&51&53&54&55&56&57&59 \\
\hline
60&61&62&63&65&66&68&69&71&73&74&75&76 \\
\hline
77&78&80&81&84&85&87&88&89&90&91&93&95 \\
\hline
96&97&98&99&100&102&103&105&106&108&109&110&112 \\
\hline
113&114&115&116&118&120&121&122&124&125&126&127&129 \\
\hline
130&131&132&133&134&135&136&138&139&141&142&143&144 \\
\hline
145&146&149&151&154&157&161&165&166&167&168&170&171 \\
\hline
175&176&177&179&180&181&182&183&184&187&189&191&192 \\
\hline
193&195&196&197&199&200&202&203&204&205&206&207&210 \\
\hline
212&214&215&216&217&220&221&223&224&225&226&227&228 \\
\hline
230&232&233&234&235&236&238&239&240&241&242&244&245 \\
\hline
248&250&251&252&253&256&257&259&262&263&264&265&266 \\
\hline
267&268&272&273&274&277&278&279&280&283&284&288&290 \\
\hline
291&292&294&295&296&297&298&299&301&304&308&309&311 \\
\hline
312&313&318&321&322&324&325&326&327&328&329&331&334 \\
\hline
335&336&337&338&339&340&341&343&344&345&347&348&349 \\
\hline
350&351&352&353&354&355&357&358&360&362&364&365&366 \\
\hline
368&369&372&373&374&378&379&380&381&383&385&386&387 \\
\hline
389&391&392&393&394&395&396&397&398&399&401&404&406 \\
\hline
407&409&410&411&412&413&414&415&416&417&418&419&420 \\
\hline
423&426&428&429&430&433&435&436&438&439&440&441&442 \\
\hline
443&446&447&451&453&454&456&457&458&460&462&465&466 \\
\hline
467&469&470&471&473&474&475&477&479&480&481&484&485 \\
\hline
486&487&489&492&496&497&498&500&501&502&504&506&507 \\
\hline
512&513&519&520&522&526&527&528&530&531&533&535&536 \\
\hline
539&540&541&543&544&549&552&553&556&558&559&561&562 \\
\hline
563&565&566&568&569&570&572&573&574&575&576&577&579 \\
\hline
580&581&582&583&590&591&593&594&596&600&601&605&606 \\
\hline
607&608&612&614&617&619&621&622&623&624&629&630&631 \\
\hline
633&634&636&640&641&644&645&646&647&648&649&650&652 \\
\hline
653&655&656&657&660&661&662&665&667&669&670&674&676 \\
\hline
678&679&680&681&684&685&687&689&693&696&697&698&699 \\
\hline
700&702&706&708&711&713&714&715&716&720&721&723&724 \\
\hline
726&728&731&733&735&736&739&742&744&747&748&749&750 \\
\hline
752&757&758&759&761&762&767&768&771&772&774&776&779 \\
\hline
780&782&783&786&787&789&791&792&793&795&797&801&802 \\
\hline
810&812&813&814&815&816&817&819&826&833&836&838&840 \\
\hline
841&845&846&851&856&860&862&864&865&867&868&870&871 \\
\hline
872&875&876&877&878&879&882&883&884&886&887 &--&--\\
\hline
\end{longtable}
\end{center}

% **********
\subsection{13a}

\begin{center}
\begin{longtable}{|c|c|c|c|c|c|c|c|c|c|c|}
\hline
3&9&13&14&17&18&19&21&23&31&34 \\
\hline
36&37&38&39&41&44&46&47&48&49&50 \\
\hline
53&54&55&60&61&64&66&67&68&69&73 \\
\hline
75&76&80&83&86&92&94&95&99&100&104 \\
\hline
107&109&110&113&114&116&120&124&126&128&139 \\
\hline
143&144&145&146&147&148&153&155&159&163&164 \\
\hline
165&169&172&174&175&176&177&178&181&183&184 \\
\hline
188&190&194&199&202&204&205&207&208&212&213 \\
\hline
216&219&222&226&227&229&233&234&235&238&239 \\
\hline
240&241&243&246&248&249&252&255&258&259&262 \\
\hline
264&266&269&270&272&276&277&280&282&283&285 \\
\hline
289&291&292&293&296&297&299&302&305&308&310 \\
\hline
312&315&320&322&323&324&325&328&329&330&332 \\
\hline
335&337&339&340&341&346&347&348&350&351&354 \\
\hline
357&359&360&363&364&367&368&371&372&373&377 \\
\hline
378&379&385&387&388&389&391&394&397&398&399 \\
\hline
400&404&409&411&413&414&415&417&418&419&420 \\
\hline
421&424&425&427&432&434&436&439&440&441&442 \\
\hline
444&446&449&450&452&453&454&456&458&459&461 \\
\hline
462&466&472&475&478&479&482&484&485&489&490 \\
\hline
491&493&499&500&502&507&509&510&511&516&517 \\
\hline
520&521&522&523&524&526&527&528&531&533&534 \\
\hline
543&547&548&549&556&557&559&560&562&563&565 \\
\hline
567&568&572&575&576&578&579&580&581&582&583 \\
\hline
586&587&588&592&595&597&600&604&605&607&608 \\
\hline
611&612&613&617&630&632&633&634&637&641&642 \\
\hline
645&646&647&649&653&659&660&666&671&672&674 \\
\hline
677&683&685&695&696&697&698&700&701&704&707 \\
\hline
710&711&715&717&719&722&723&727&728&733&735 \\
\hline
736&739&741&743&745&746&748&749&750&751&759 \\
\hline
760&761&764&766&767&769&770&771&778&780&781 \\
\hline
782&785&788&790&791&792&793&794&795&796&798 \\
\hline
799&800&801&802&803&807&809&811&812&817&820 \\
\hline
821&822&825&828&829&833&837&842&843&844&847 \\
\hline
848&849&852&853&856&858&859&860&862&863&868 \\
\hline
869&872&875&876&878&882&884&886&887&888&891 \\
\hline
893&898&899&904&907&908&912&914&915&917&920 \\
\hline
921&929&930&932&934&935&937&939&942&944&946 \\
\hline
948&949&953&954&958&959&965&969&972&975&977 \\
\hline
978&979&983&984&987&992&994&996&997&998&999 \\
\hline
1001&1005&1007&1009&1010&1012&1014&1016&1017&1023&1024 \\
\hline
1025&1026&1029&1030&1031&1034&1036&1037&1039&1040&1042 \\
\hline
1044&1046&1047&1049&1052&1053&1054&1055&1056&1057&1061 \\
\hline
1062&1064&1065&1071&1073&1074&1079&1080&1081&1082&1083 \\
\hline
1084&1086&1088&1094&1095&1096&1097&1098&1101&1105&1108 \\
\hline
1109&1113&1114&1115&1116&1117&1119&1124&1125&1126&1128 \\
\hline
1130&1132&1136&1138&1139&1141&1143&1144&1146&1148&1150 \\
\hline
1154&1156&1160&1162&1164&1165&1166&1167&1168&1170&1171 \\
\hline
1172&1174&1175&1177&1178&1181&1184&1185&1187&1189&1190 \\
\hline
1191&1194&1195&1199&1201&1202&1204&1205&1206&1208&1209 \\
\hline
1211&1215&1218&1220&1223&1224&1225&1227&1228&1229&1236 \\
\hline
1241&1244&1248&1249&1252&1253&1254&1256&1257&1260&1262 \\
\hline
1263&1268&1270&1273&1275&1276&1278&1280&1282&1284&1286 \\
\hline
1287&1288&1289&1290&1291&1294&1295&1296&1297&1300&1303 \\
\hline
1304&1308&1309&1314&1315&1317&1318&1319&1322&1323&1327 \\
\hline
1329&1333&1334&1336&1337&1338&1340&1341&1343&1344&1347 \\
\hline
1352&1354&1356&1358&1360&1361&1363&1365&1366&1370&1373 \\
\hline
1374&1375&1376&1379&1382&1383&1386&1387&1388&1389&1390 \\
\hline
1391&1392&1393&1397&1399&1401&1402&1403&1407&1409&1411 \\
\hline
1413&1415&1416&1420&1423&1424&1426&1427&1432&1433&1435 \\
\hline
1441&1445&1449&1450&1451&1453&1458&1459&1469&1470&1473 \\
\hline
1476&1478&1479&1480&1484&1487&1489&1491&1493&1494&1497 \\
\hline
1498&1499&1501&1502&1505&1507&1508&1509&1512&1513&1516 \\
\hline
1518&1519&1521&1527&1529&1530&1534&1536&1537&1538&1540 \\
\hline
1545&1549&1551&1552&1554&1555&1558&1559&1560&1561&1562 \\
\hline
1563&1564&1568&1572&1575&1576&1577&1579&1581&1582&1584 \\
\hline
1585&1586&1587&1588&1591&1595&1599&1600&1605&1606&1607 \\
\hline
1611&1612&1613&1615&1616&1619&1620&1628&1631&1635&1636 \\
\hline
1639&1642&1644&1646&1647&1650&1653&1657&1663&1664&1669 \\
\hline
1670&1674&1676&1678&1679&1680&1681&1684&1685&1687&1688 \\
\hline
1690&1691&1693&1699&1700&1702&1705&1707&1708&1710&1711 \\
\hline
1713&1715&1716&1726&1729&1730&1733&1734&1736&1737&1738 \\
\hline
1739&1741&1742&1744&1745&1746&1747&1749&1752&1755&1759 \\
\hline
1760&1764&1765&1769&1770&1771&1772&1773&1774&1776&1778 \\
\hline
1780&1782&1784&1788&1790&1793&1795&1796&1798&1799&1800 \\
\hline
1801&1802&1803&1807&1808&1809&1810&1812&1813&1815&1818 \\
\hline
1821&1822&1823&1824&1828&1831&1832&1834&1836&1838&1840 \\
\hline
1842&1844&1846&1847&1848&1850&1853&1858&1860&1861&1868 \\
\hline
1870&1874&1876&1878&1880&1881&1882&1883&1884&1886&1887 \\
\hline
1888&1890&1891&1892&1893&1898&1900&1902&1903&1905&1910 \\
\hline
1911&1913&1914&1915&1916&1919&1923&1924&1925&1926&1930 \\
\hline
1931&1932&1933&1938&1939&1940&1941&1942&1944&1948&1950 \\
\hline
1953&1954&1958&1959&1964&1965&1967&1968&1973&1977&1980 \\
\hline
1981&1983&1994&1995&1996&1998&2000&2001&2002&2003&2004 \\
\hline
2005&2006&2009&2012&2016&2017&2022&2023&2033&2034&2035 \\
\hline
2039&2041&2043&2044&2045&2047&2049&2050&2051&2053&2056 \\
\hline
2057&2062&2065&2066&2067&2068&2069&2070&2072&2084&2085 \\
\hline
2086&2087&2089&2090&2095&2099&2100&2102&2105&2110&2111 \\
\hline
2115&2117&2118&2120&2121&2122&2123&2125&2128&2134&2136 \\
\hline
2137&2140&2142&2143&2144&2145&2146&2149&2151&2152&2155 \\
\hline
2158&2160&2161&2162&2165&2168&2169&2172&2174&2178&2179 \\
\hline
2182&2183&2184&2186&2187&2190&2192&2195&2197&2198&2203 \\
\hline
2206&2208&2209&2211&2212&2215&2219&2222&2223&2226&2227 \\
\hline
2229&2230&2231&2235&2236&2239&2240&2241&2242&2243&2246 \\
\hline
2248&2250&2251&2252&2253&2256&2258&2259&2260&2262&2264 \\
\hline
2267&2272&2274&2275&2276&2277&2279&2280&2283&2285&2288 \\
\hline
2289&2291&2293&2294&2295&2297&2298&2299&2300&2302&2303 \\
\hline
2304&2309&2312&2313&2314&2316&2317&2320&2321&2325&2327 \\
\hline
2329&2336&2338&2340&2342&2343&2344&2345&2347&2348&2349 \\
\hline
2350&2353&2354&2355&2356&2357&2358&2361&2363&2364&2365 \\
\hline
2368&2369&2371&2373&2377&2379&2380&2382&2383&2384&2385 \\
\hline
2391&2393&2396&2397&2401&2403&2404&2410&2412&2413&2415 \\
\hline
2416&2418&2422&2423&2426&2427&2429&2430&2431&2435&2437 \\
\hline
2439&2440&2441&2444&2446&2447&2448&2449&2451&2454&2456 \\
\hline
2459&2463&2464&2465&2466&2468&2472&2473&2475&2477&2479 \\
\hline
2480&2482&2483&2486&2490&2492&2494&2496&2498&2504&2508 \\
\hline
2510&2513&2515&2516&2517&2518&2519&2521&2522&2527&2528 \\
\hline
2529&2532&2537&2539&2540&2545&2546&2547&2549&2551&2553 \\
\hline
2558&2559&2566&2567&2568&2569&2570&2572&2574&2576&2579 \\
\hline
2585&2588&2596&2599&2602&2603&2609&2612&2613&2614&2617 \\
\hline
2621&2623&2625&2634&2638&2648&2661&2664&2665&2675&2677 \\
\hline
2680&2681&2688&2689&2692&2693&2694&2696&2698&2700&2702 \\
\hline
2704&2705&2708&2712&2713&2714&2725&2728&2729&2730&2731 \\
\hline
2733&2734&2740&2741&2743&2746&2749&2754&2755&2756&2758 \\
\hline
2759&2760&2764&2766&2768&2770&2774&2775&2776&2777&2779 \\
\hline
2780&2783&2785&2789&2793&2794&2795&2799&2803&2804&2805 \\
\hline
2807&2809&2812&2817&2818&2819&2820&2822&2825&2827&2831 \\
\hline
2832&2834&2836&2837&2839&2841&2842&2845&2847&2848&2850 \\
\hline
2851&2852&2854&2857&2858&2859&2861&2862&2867&2868&2869 \\
\hline
2871&2874&2882&2883&2885&2891&2894&2899&2902&2903&2904 \\
\hline
2905&2906&2908&2909&2910&2912&2913&2914&2917&2921&2922 \\
\hline
2926&2927&2928&2929&2931&2932&2934&2938&2941&2942&2948 \\
\hline
2952&2954&2955&2956&2957&2961&2963&2964&2965&2973&2974 \\
\hline
2975&2976&2977&2978&2979&2980&2986&2989&2990&2991&2992 \\
\hline
2993&2997&2998&3002&3003&3004&3006&3019&3023&3029&3031 \\
\hline
3032&3033&3035&3038&3041&3047&3048&3052&3053&3054&3064 \\
\hline
3065&3067&3069&3070&3071&3073&3077&3078&3081&3082&3084 \\
\hline
3085&3086&3092&3094&3095&3097&3098&3101&3102&3103&3105 \\
\hline
3106&3109&3111&3113&3115&3117&3119&3123&3124&3134&3136 \\
\hline
3137&3154&3162&3163&3165&3169&3171&3175&3183&3184&3186 \\
\hline
3189&3193&3197&3201&3202&3203&3205&3208&3209&3211&3213 \\
\hline
3214&3215&3218&3220&3221&3227&3228&3229&3233&3241&3243 \\
\hline
3246&3247&3251&3252&3260&3264&3269&3270&3271&3276&3278 \\
\hline
3280&3286&3287&3293&3295&3299&3305&3310&3311&3312&3314 \\
\hline
3316&3317&3318&3319&3320&3321&3322&3325&3330&3331&3335 \\
\hline
3338&3342&3345&3346&3347&3348&3350&3357&3359&3363&3366 \\
\hline
3368&3370&3378&3381&3382&3384&3385&3389&3391&3394&3398 \\
\hline
3407&3411&3412&3416&3420&3421&3423&3427&3429&3430&3431 \\
\hline
3432&3433&3434&3437&3439&3440&3443&3447&3448&3450&3451 \\
\hline
3452&3453&3456&3464&3469&3470&3471&3473&3480&3481&3483 \\
\hline
3484&3488&3496&3498&3499&3503&3507&3513&3518&3519&3523 \\
\hline
3525&3526&3527&3532&3539&3540&3541&3545&3548&3549&3550 \\
\hline
3555&3560&3561&3562&3565&3569&3574&3575&3576&3580&3584 \\
\hline
3587&3593&3599&3603&3606&3607&3610&3614&3616&3617&3619 \\
\hline
3625&3628&3629&3630&3637&3640&3643&3644&3646&3647&3649 \\
\hline
3654&3656&3661&3662&3664&3666&3668&3671&3673&3675&3679 \\
\hline
3680&3681&3682&3694&3695&3697&3698&3700&3703&3705&3713 \\
\hline
3715&3719&3722&3723&3724&3725&3727&3728&3729&3732&3733 \\
\hline
3738&3739&3744&3755&3756&3760&3761&3766&3767&3771&3773 \\
\hline
3774&3778&3779&3780&3781&3783&3784&3788&3789&3790&3791 \\
\hline
3797&3802&3806&3808&3813&3814&3815&3816&3817&3818&3820 \\
\hline
3821&3824&3825&3829&3830&3831&3832&3833&3834&3835&3838 \\
\hline
3849&3851&3856&3859&3860&3861&3863&3867&3870&3871&3874 \\
\hline
3876&3878&3881&3886&3889&3893&3894&3895&3897&3898&3900 \\
\hline
3902&3903&3908&3909&3910&3915&3918&3919&3923&3924&3929 \\
\hline
3932&3933&3934&3935&3938&3949&3952&3963&3964&3965&3968 \\
\hline
3971&3974&3980&3981&3984&3985&3987&3988&3991&3994&3997 \\
\hline
4007&4008&4013&4014&4016&4020&4022&4023&4024&4026&4027 \\
\hline
4029&4033&4036&4040&4042&4061&4064&4067&4068&4071&4074 \\
\hline
4075&4076&4078&4081&4082&4084&4085&4086&4091&4094&4095 \\
\hline
4105&4111&4114&4115&4117&4121&4122&4124&4131&4133&4135 \\
\hline
4137&4138&4139&4142&4145&4151&4156&4159&4161&4164&4165 \\
\hline
4166&4169&4171&4175&4180&4181&4185&4186&4187&4189&4196 \\
\hline
4197&4199&4200&4202&4203&4204&4213&4215&4221&4222&4223 \\
\hline
4227&4230&4232&4236&4237&4242&4243&4244&4245&4256&4261 \\
\hline
4262&4264&4271&4273&4275&4276&4277&4278&4279&4280&4288 \\
\hline
4291&4292&4294&4295&4296&4297&4299&4300&4302&4303&4305 \\
\hline
4311&4314&4315&4317&4318&4319&4323&4324&4326&4332&4335 \\
\hline
4336&4340&4347&4351&4354&4357&4364&4370&4372&4373&4375 \\
\hline
4378&4383&4384&4387&4395&4396&4399&4400&4402&4405&4407 \\
\hline
4411&4413&4415&4417&4425&4428&4432&4433&4434&4439&4440 \\
\hline
4446&4447&4450&4453&4455&4458&4461&4462&4469&4470&4477 \\
\hline
4485&4489&4491&4496&4500&4501&4510&4515&4520&4521&4522 \\
\hline
4527&4533&4534&4539&4542&4543&4547&4548&4550&4556&4560 \\
\hline
4561&4562&4563&4564&4569&4570&4572&4573&4579&4590&4593 \\
\hline
4596&4597&4598&4602&4605&4608&4610&4613&4617&4618&4619 \\
\hline
4623&4627&4628&4630&4633&4634&4635&4638&4639&4640&4641 \\
\hline
4642&4651&4652&4655&4659&4661&4662&4664&4666&4668&4675 \\
\hline
4676&4679&4682&4687&4689&4699&4700&4702&4703&4708&4710 \\
\hline
4711&4721&4723&4726&4728&4734&4738&4741&4748&4750&4751 \\
\hline
4752&4753&4759&4760&4772&4773&4775&4776&4779&4780&4781 \\
\hline
4783&4784&4785&4788&4789&4793&4795&4797&4801&4802&4803 \\
\hline
4804&4805&4808&4810&4814&4818&4819&4822&4825&4827&4829 \\
\hline
4830&4831&4833&4834&4836&4839&4841&4852&4854&4855&4863 \\
\hline
4866&4867&4868&4870&4871&4872&4878&--  &--  &--  &-- \\
\hline
\end{longtable}
\end{center}

% **********
\subsection{13n}

\begin{center}
\begin{longtable}{|c|c|c|c|c|c|c|c|c|c|c|}
\hline
1&3&4&5&10&11&13&14&15&18&21 \\
\hline
22&24&26&27&28&29&31&33&35&36&37 \\
\hline
38&40&43&46&47&48&49&51&52&53&55 \\
\hline
56&57&59&61&62&65&67&68&69&71&72 \\
\hline
74&76&77&79&80&81&83&84&86&88&91 \\
\hline
94&96&98&100&102&106&107&108&109&110&111 \\
\hline
112&114&116&117&118&120&121&122&123&124&125 \\
\hline
126&128&130&132&133&134&135&136&138&140&141 \\
\hline
143&144&145&146&147&148&151&155&156&158&160 \\
\hline
161&162&164&165&167&168&169&170&171&172&173 \\
\hline
174&175&176&179&180&181&182&183&184&185&188 \\
\hline
189&190&191&192&194&195&197&198&200&202&203 \\
\hline
206&207&208&209&211&212&213&214&216&220&221 \\
\hline
223&224&226&228&230&232&233&234&235&236&237 \\
\hline
239&240&244&246&251&253&255&257&258&259&260 \\
\hline
261&262&263&264&267&268&269&270&273&275&279 \\
\hline
280&281&282&283&284&287&288&291&292&293&294 \\
\hline
296&297&299&303&305&309&310&311&314&316&318 \\
\hline
319&320&322&323&325&327&328&329&332&333&335 \\
\hline
337&340&343&345&349&350&351&354&355&357&358 \\
\hline
359&361&362&363&365&367&369&371&372&373&374 \\
\hline
375&379&381&382&387&389&391&392&395&398&399 \\
\hline
400&402&403&405&407&408&409&410&411&412&413 \\
\hline
414&415&417&419&420&421&422&423&424&425&430 \\
\hline
431&433&434&435&436&440&442&445&447&449&450 \\
\hline
452&454&455&457&459&463&464&467&469&470&471 \\
\hline
473&475&477&478&481&482&484&485&487&492&494 \\
\hline
496&497&499&502&504&505&506&509&511&512&513 \\
\hline
516&518&521&523&525&526&528&531&532&534&536 \\
\hline
539&540&542&543&544&545&547&548&549&551&552 \\
\hline
553&555&557&558&562&564&566&567&572&577&578 \\
\hline
579&581&583&584&586&587&588&589&592&593&599 \\
\hline
602&603&604&607&608&609&611&612&613&621&623 \\
\hline
624&626&630&632&636&637&640&641&642&643&644 \\
\hline
647&648&649&653&654&655&656&657&658&659&660 \\
\hline
664&665&668&669&673&676&677&678&682&683&684 \\
\hline
685&687&690&692&696&699&700&701&703&706&707 \\
\hline
708&711&717&722&724&725&726&727&729&730&731 \\
\hline
735&737&739&740&741&743&745&749&750&751&756 \\
\hline
757&758&760&762&763&765&766&768&770&775&777 \\
\hline
779&781&783&784&785&786&788&791&792&793&794 \\
\hline
795&797&800&801&802&806&807&809&810&814&815 \\
\hline
818&820&823&824&826&827&833&835&836&840&842 \\
\hline
844&848&851&853&856&857&859&861&862&863&866 \\
\hline
867&868&871&872&876&877&879&880&881&882&883 \\
\hline
886&887&888&890&894&896&897&898&901&904&905 \\
\hline
908&910&911&912&913&915&917&918&919&920&921 \\
\hline
922&924&928&930&931&933&936&937&939&940&941 \\
\hline
944&945&948&951&952&953&954&955&957&958&959 \\
\hline
963&964&967&971&973&976&977&979&980&982&986 \\
\hline
987&989&993&995&998&1001&1007&1010&1012&1014&1015 \\
\hline
1016&1017&1019&1021&1022&1023&1025&1026&1027&1028&1030 \\
\hline
1031&1032&1034&1038&1040&1044&1046&1048&1051&1052&1053 \\
\hline
1054&1056&1059&1060&1064&1065&1068&1070&1071&1073&1074 \\
\hline
1076&1081&1082&1084&1087&1088&1089&1090&1091&1093&1096 \\
\hline
1098&1100&1101&1102&1103&1105&1107&1108&1109&1110&1112 \\
\hline
1116&1118&1120&1121&1122&1125&1127&1132&1135&1136&1139 \\
\hline
1142&1144&1146&1147&1149&1150&1151&1152&1153&1155&1156 \\
\hline
1158&1159&1160&1166&1167&1168&1169&1170&1171&1172&1173 \\
\hline
1174&1175&1178&1181&1182&1186&1189&1190&1191&1194&1195 \\
\hline
1196&1198&1199&1200&1206&1207&1209&1214&1216&1218&1219 \\
\hline
1221&1223&1224&1226&1228&1230&1232&1233&1237&1241&1245 \\
\hline
1246&1247&1248&1251&1252&1254&1256&1257&1260&1261&1262 \\
\hline
1264&1265&1266&1267&1268&1269&1271&1274&1276&1280&1282 \\
\hline
1283&1285&1286&1287&1291&1292&1293&1295&1296&1298&1299 \\
\hline
1301&1303&1304&1305&1306&1308&1309&1310&1312&1313&1314 \\
\hline
1316&1317&1319&1320&1321&1322&1324&1325&1326&1327&1331 \\
\hline
1334&1336&1337&1338&1340&1342&1343&1345&1349&1351&1352 \\
\hline
1353&1357&1359&1360&1363&1366&1367&1368&1371&1372&1373 \\
\hline
1376&1378&1381&1384&1386&1387&1388&1391&1392&1393&1395 \\
\hline
1398&1400&1404&1405&1407&1410&1412&1413&1416&1418&1420 \\
\hline
1423&1424&1425&1426&1427&1428&1429&1430&1431&1432&1433 \\
\hline
1434&1435&1436&1442&1443&1444&1445&1447&1448&1449&1450 \\
\hline
1451&1452&1453&1455&1456&1457&1458&1459&1464&1465&1466 \\
\hline
1468&1470&1471&1472&1474&1475&1478&1479&1480&1481&1482 \\
\hline
1485&1487&1491&1493&1495&1496&1499&1501&1505&1510&1512 \\
\hline
1513&1514&1521&1522&1527&1531&1532&1534&1536&1538&1540 \\
\hline
1541&1544&1545&1546&1547&1548&1549&1550&1556&1559&1561 \\
\hline
1562&1563&1564&1565&1566&1569&1571&1573&1574&1575&1578 \\
\hline
1580&1583&1587&1588&1593&1596&1597&1598&1599&1600&1601 \\
\hline
1603&1604&1605&1607&1608&1611&1612&1613&1614&1615&1616 \\
\hline
1617&1619&1620&1622&1623&1624&1625&1626&1627&1630&1631 \\
\hline
1632&1634&1635&1636&1637&1641&1642&1644&1645&1646&1649 \\
\hline
1650&1651&1652&1653&1655&1658&1659&1661&1662&1663&1664 \\
\hline
1666&1667&1668&1671&1674&1675&1680&1682&1683&1685&1690 \\
\hline
1691&1692&1694&1695&1696&1697&1698&1699&1701&1702&1706 \\
\hline
1709&1710&1711&1712&1713&1715&1717&1719&1721&1722&1723 \\
\hline
1724&1725&1727&1728&1729&1731&1733&1734&1735&1737&1738 \\
\hline
1740&1741&1742&1744&1746&1749&1752&1756&1757&1759&1760 \\
\hline
1763&1764&1765&1766&1767&1769&1772&1775&1777&1779&1780 \\
\hline
1788&1792&1793&1794&1795&1797&1798&1803&1805&1806&1807 \\
\hline
1808&1810&1811&1812&1813&1814&1815&1819&1821&1823&1824 \\
\hline
1825&1826&1829&1830&1831&1834&1835&1836&1837&1841&1845 \\
\hline
1846&1847&1848&1849&1852&1853&1854&1856&1857&1859&1860 \\
\hline
1865&1867&1868&1869&1871&1874&1875&1877&1878&1879&1880 \\
\hline
1881&1882&1883&1884&1885&1886&1887&1888&1889&1890&1891 \\
\hline
1893&1896&1897&1898&1899&1900&1901&1902&1904&1905&1911 \\
\hline
1912&1914&1916&1919&1921&1922&1924&1926&1928&1929&1930 \\
\hline
1934&1937&1940&1941&1942&1943&1945&1948&1952&1953&1954 \\
\hline
1955&1956&1958&1959&1960&1961&1963&1968&1970&1971&1972 \\
\hline
1973&1974&1975&1979&1981&1983&1984&1985&1986&1988&1989 \\
\hline
1990&1991&1993&1995&1996&1997&1999&2000&2001&2004&2005 \\
\hline
2006&2007&2008&2011&2012&2014&2015&2016&2017&2019&2020 \\
\hline
2023&2024&2027&2028&2029&2030&2031&2033&2034&2036&2037 \\
\hline
2038&2039&2040&2043&2046&2047&2048&2049&2051&2052&2054 \\
\hline
2055&2058&2060&2062&2064&2066&2068&2069&2070&2072&2073 \\
\hline
2075&2076&2078&2079&2081&2082&2084&2087&2090&2091&2092 \\
\hline
2093&2094&2095&2097&2098&2100&2101&2102&2105&2106&2112 \\
\hline
2113&2120&2122&2124&2126&2127&2128&2129&2130&2132&2133 \\
\hline
2138&2140&2142&2144&2145&2148&2149&2151&2152&2154&2156 \\
\hline
2157&2158&2161&2162&2163&2164&2165&2166&2167&2168&2169 \\
\hline
2172&2173&2175&2176&2177&2178&2179&2180&2181&2185&2187 \\
\hline
2188&2189&2194&2196&2197&2198&2199&2200&2201&2204&2205 \\
\hline
2209&2210&2212&2214&2215&2216&2217&2219&2220&2223&2224 \\
\hline
2227&2228&2229&2233&2234&2235&2236&2238&2239&2240&2241 \\
\hline
2242&2243&2244&2246&2247&2250&2254&2255&2256&2257&2259 \\
\hline
2260&2261&2264&2266&2268&2269&2273&2274&2275&2277&2281 \\
\hline
2282&2283&2285&2286&2289&2290&2293&2295&2298&2299&2300 \\
\hline
2301&2303&2306&2307&2308&2309&2310&2315&2316&2317&2318 \\
\hline
2322&2323&2324&2325&2328&2333&2336&2339&2342&2344&2347 \\
\hline
2348&2353&2354&2355&2356&2357&2358&2359&2360&2363&2364 \\
\hline
2367&2370&2371&2374&2376&2378&2379&2380&2382&2384&2385 \\
\hline
2387&2388&2389&2390&2392&2393&2394&2396&2398&2400&2401 \\
\hline
2407&2410&2412&2414&2416&2417&2420&2422&2423&2424&2427 \\
\hline
2428&2429&2430&2431&2432&2435&2436&2441&2443&2445&2446 \\
\hline
2447&2448&2449&2451&2452&2456&2457&2459&2462&2463&2464 \\
\hline
2466&2467&2470&2474&2475&2476&2477&2479&2480&2481&2491 \\
\hline
2492&2493&2495&2496&2497&2500&2502&2504&2506&2507&2508 \\
\hline
2509&2511&2512&2513&2514&2517&2518&2519&2520&2522&2523 \\
\hline
2524&2525&2531&2533&2534&2537&2538&2540&2542&2543&2544 \\
\hline
2545&2546&2549&2552&2553&2554&2556&2557&2558&2559&2560 \\
\hline
2562&2566&2567&2568&2569&2570&2571&2573&2577&2578&2579 \\
\hline
2580&2581&2582&2584&2589&2590&2592&2593&2594&2595&2597 \\
\hline
2598&2600&2601&2602&2604&2605&2613&2614&2616&2617&2619 \\
\hline
2620&2621&2623&2624&2625&2630&2631&2634&2637&2638&2640 \\
\hline
2644&2647&2648&2650&2651&2652&2655&2656&2660&2662&2663 \\
\hline
2666&2667&2668&2669&2673&2674&2678&2679&2680&2681&2682 \\
\hline
2684&2685&2686&2687&2691&2692&2695&2696&2697&2698&2700 \\
\hline
2701&2702&2703&2704&2705&2707&2708&2709&2710&2712&2713 \\
\hline
2714&2715&2716&2717&2718&2720&2721&2723&2725&2726&2727 \\
\hline
2729&2731&2734&2735&2737&2738&2739&2742&2744&2747&2748 \\
\hline
2749&2752&2753&2754&2756&2759&2762&2768&2769&2771&2772 \\
\hline
2774&2776&2778&2779&2780&2781&2782&2785&2786&2787&2791 \\
\hline
2792&2795&2796&2797&2799&2801&2802&2803&2806&2807&2809 \\
\hline
2810&2813&2815&2817&2818&2820&2823&2824&2825&2826&2831 \\
\hline
2832&2833&2835&2836&2841&2842&2843&2844&2846&2847&2848 \\
\hline
2849&2852&2853&2856&2859&2864&2865&2866&2867&2868&2869 \\
\hline
2872&2876&2878&2880&2882&2884&2885&2886&2887&2888&2889 \\
\hline
2890&2894&2896&2897&2902&2903&2904&2906&2908&2910&2911 \\
\hline
2912&2913&2915&2917&2918&2919&2921&2922&2924&2927&2928 \\
\hline
2930&2933&2934&2937&2940&2941&2943&2944&2947&2950&2951 \\
\hline
2953&2954&2955&2956&2957&2958&2960&2962&2964&2965&2969 \\
\hline
2970&2971&2973&2974&2975&2978&2984&2987&2988&2989&2990 \\
\hline
2991&2992&2993&2996&2997&2999&3002&3006&3007&3008&3009 \\
\hline
3011&3012&3013&3014&3016&3019&3021&3022&3024&3026&3027 \\
\hline
3028&3029&3030&3031&3032&3033&3035&3036&3037&3039&3040 \\
\hline
3041&3042&3044&3047&3050&3052&3053&3054&3055&3058&3061 \\
\hline
3062&3064&3065&3066&3067&3069&3070&3071&3074&3075&3076 \\
\hline
3077&3078&3079&3081&3082&3084&3085&3086&3087&3090&3091 \\
\hline
3092&3094&3097&3098&3099&3100&3101&3105&3110&3112&3114 \\
\hline
3115&3118&3119&3120&3125&3126&3127&3128&3132&3134&3135 \\
\hline
3138&3142&3143&3144&3148&3150&3151&3152&3153&3157&3158 \\
\hline
3162&3163&3165&3167&3169&3171&3172&3173&3175&3176&3177 \\
\hline
3179&3180&3184&3185&3186&3187&3191&3194&3195&3196&3197 \\
\hline
3201&3203&3204&3205&3206&3208&3211&3214&3220&3221&3222 \\
\hline
3225&3227&3228&3229&3231&3235&3238&3240&3241&3242&3243 \\
\hline
3244&3246&3247&3250&3251&3256&3257&3259&3260&3261&3264 \\
\hline
3267&3268&3270&3272&3275&3279&3280&3281&3284&3289&3290 \\
\hline
3291&3295&3297&3300&3301&3302&3306&3308&3312&3313&3314 \\
\hline
3316&3317&3318&3319&3320&3321&3322&3324&3326&3329&3332 \\
\hline
3335&3340&3342&3343&3345&3346&3347&3348&3350&3355&3356 \\
\hline
3359&3361&3362&3363&3364&3365&3370&3371&3372&3374&3375 \\
\hline
3376&3378&3379&3380&3383&3386&3389&3390&3391&3392&3393 \\
\hline
3394&3395&3396&3400&3401&3403&3404&3405&3406&3408&3409 \\
\hline
3410&3411&3413&3414&3415&3416&3417&3418&3420&3421&3425 \\
\hline
3429&3431&3433&3436&3438&3443&3444&3448&3450&3451&3453 \\
\hline
3456&3457&3462&3463&3464&3465&3466&3467&3468&3469&3471 \\
\hline
3472&3473&3474&3477&3478&3480&3482&3487&3488&3490&3494 \\
\hline
3495&3500&3503&3504&3506&3509&3510&3511&3513&3516&3517 \\
\hline
3519&3522&3523&3524&3525&3527&3528&3531&3533&3535&3537 \\
\hline
3539&3541&3542&3543&3545&3549&3551&3552&3554&3555&3557 \\
\hline
3559&3560&3565&3567&3569&3571&3574&3575&3577&3578&3580 \\
\hline
3582&3585&3586&3588&3590&3592&3594&3598&3601&3604&3605 \\
\hline
3606&3607&3609&3610&3612&3614&3617&3619&3621&3622&3624 \\
\hline
3625&3629&3630&3631&3632&3634&3635&3638&3642&3644&3646 \\
\hline
3649&3650&3651&3654&3655&3665&3666&3668&3669&3670&3671 \\
\hline
3675&3677&3678&3679&3680&3681&3682&3683&3684&3685&3686 \\
\hline
3687&3688&3689&3690&3693&3694&3696&3700&3701&3702&3703 \\
\hline
3705&3706&3707&3708&3711&3713&3715&3716&3717&3725&3726 \\
\hline
3727&3728&3730&3734&3737&3738&3741&3742&3746&3747&3748 \\
\hline
3753&3754&3756&3758&3759&3760&3761&3763&3765&3766&3769 \\
\hline
3770&3772&3773&3777&3778&3779&3783&3784&3787&3789&3792 \\
\hline
3793&3796&3797&3798&3800&3802&3803&3809&3811&3812&3814 \\
\hline
3815&3818&3819&3820&3821&3823&3825&3828&3830&3834&3835 \\
\hline
3838&3839&3841&3842&3845&3846&3847&3850&3851&3853&3856 \\
\hline
3858&3859&3861&3862&3863&3865&3866&3867&3868&3869&3870 \\
\hline
3871&3872&3874&3875&3876&3880&3881&3883&3885&3887&3888 \\
\hline
3889&3894&3895&3897&3900&3904&3905&3908&3909&3910&3911 \\
\hline
3912&3914&3915&3917&3925&3926&3927&3930&3931&3932&3934 \\
\hline
3936&3937&3938&3939&3940&3943&3946&3947&3948&3949&3950 \\
\hline
3952&3953&3954&3955&3956&3957&3960&3962&3963&3964&3965 \\
\hline
3967&3971&3972&3974&3975&3977&3979&3980&3982&3983&3985 \\
\hline
3986&3988&3991&3992&3993&3994&3996&3997&3998&4001&4002 \\
\hline
4003&4004&4005&4007&4009&4010&4011&4012&4014&4016&4017 \\
\hline
4018&4019&4022&4024&4026&4027&4028&4031&4032&4033&4034 \\
\hline
4036&4037&4038&4039&4040&4041&4045&4046&4049&4051&4052 \\
\hline
4053&4056&4058&4061&4063&4064&4067&4074&4075&4076&4079 \\
\hline
4080&4081&4082&4083&4084&4086&4088&4089&4090&4092&4096 \\
\hline
4097&4098&4099&4101&4102&4103&4105&4106&4108&4109&4112 \\
\hline
4114&4116&4118&4120&4121&4122&4124&4126&4127&4131&4136 \\
\hline
4138&4140&4142&4143&4144&4145&4146&4147&4149&4150&4151 \\
\hline
4157&4160&4161&4163&4166&4168&4170&4172&4174&4175&4178 \\
\hline
4184&4185&4188&4190&4191&4194&4196&4198&4201&4203&4204 \\
\hline
4206&4209&4210&4212&4214&4220&4222&4226&4227&4228&4230 \\
\hline
4231&4232&4233&4235&4236&4237&4242&4244&4248&4249&4251 \\
\hline
4253&4254&4255&4257&4258&4259&4262&4264&4265&4266&4270 \\
\hline
4271&4273&4274&4275&4276&4279&4281&4282&4285&4286&4287 \\
\hline
4288&4293&4295&4296&4297&4298&4299&4302&4303&4304&4311 \\
\hline
4312&4320&4321&4322&4324&4326&4327&4328&4331&4332&4335 \\
\hline
4336&4340&4347&4348&4352&4355&4356&4357&4364&4368&4374 \\
\hline
4376&4378&4382&4384&4386&4390&4392&4393&4397&4399&4400 \\
\hline
4401&4403&4405&4406&4407&4409&4411&4412&4419&4421&4423 \\
\hline
4424&4425&4427&4428&4429&4430&4432&4436&4439&4440&4443 \\
\hline
4446&4453&4456&4461&4462&4464&4465&4466&4470&4471&4472 \\
\hline
4473&4474&4476&4480&4481&4482&4484&4486&4487&4488&4489 \\
\hline
4492&4494&4498&4500&4502&4503&4505&4508&4510&4512&4513 \\
\hline
4514&4516&4519&4523&4524&4525&4527&4528&4531&4533&4534 \\
\hline
4535&4538&4540&4541&4542&4545&4548&4551&4552&4558&4559 \\
\hline
4560&4562&4563&4566&4573&4574&4578&4579&4580&4581&4582 \\
\hline
4583&4584&4585&4587&4590&4591&4592&4593&4594&4596&4598 \\
\hline
4602&4608&4610&4612&4613&4617&4621&4623&4624&4625&4626 \\
\hline
4628&4629&4631&4634&4635&4639&4640&4643&4646&4648&4654 \\
\hline
4655&4658&4659&4662&4664&4667&4669&4670&4671&4675&4677 \\
\hline
4681&4682&4683&4684&4685&4690&4691&4695&4699&4701&4702 \\
\hline
4705&4706&4710&4712&4713&4715&4716&4718&4720&4721&4722 \\
\hline
4723&4724&4725&4727&4732&4735&4736&4738&4740&4741&4742 \\
\hline
4755&4758&4759&4763&4766&4768&4769&4771&4772&4776&4780 \\
\hline
4783&4787&4788&4789&4790&4792&4794&4796&4797&4800&4801 \\
\hline
4803&4805&4806&4807&4808&4809&4816&4819&4821&4824&4825 \\
\hline
4826&4827&4829&4831&4832&4833&4834&4835&4836&4839&4840 \\
\hline
4843&4844&4845&4849&4850&4851&4854&4855&4858&4859&4862 \\
\hline
4864&4869&4871&4874&4875&4876&4877&4881&4882&4883&4885 \\
\hline
4893&4894&4895&4901&4902&4905&4907&4908&4909&4912&4915 \\
\hline
4917&4918&4920&4921&4923&4924&4929&4932&4934&4937&4939 \\
\hline
4944&4945&4946&4948&4950&4954&4956&4957&4959&4961&4962 \\
\hline
4965&4968&4969&4973&4976&4979&4986&4987&4988&4992&4995 \\
\hline
4996&5001&5002&5004&5007&5011&5016&5017&5018&5019&5020 \\
\hline
5023&5026&5028&5029&5031&5032&5033&5034&5048&5049&5050 \\
\hline
5051&5053&5056&5057&5058&5061&5063&5064&5065&5069&5071 \\
\hline
5074&5080&5081&5082&5086&5088&5090&5091&5092&5093&5094 \\
\hline
5096&5099&5103&5105&5106&5110&--  &--  &--  &--  &-- \\
\hline
\end{longtable}
\end{center}

% ********************
\bibliographystyle{amsplain}
\bibliography{main}

\end{document}